\documentclass{article}

\usepackage{amsthm}
\usepackage{amsmath}
\usepackage{amssymb}
\usepackage{amscd}                
\usepackage[all]{xy}               
\usepackage{stmaryrd}              
\usepackage{a4wide}

\theoremstyle{plain}

\newtheorem{theorem}{Theorem}[section]
\newtheorem{lemma}[theorem]{Lemma}

\theoremstyle{definition}

\newtheorem{definition}[theorem]{Definition}
\newtheorem{example}[theorem]{Example}

\newtheorem{prob}[theorem]{Problem}
\newtheorem{setup}[theorem]{Set-up}

\numberwithin{equation}{section}

\newcommand\C{{\mathbb C}}
\newcommand\Q{{\mathbb Q}}
\newcommand\Z{{\mathbb Z}}

\newcommand{\Gal}{\mathop{\mathrm{Gal} }\nolimits}
\newcommand{\GL}{\mathop{\mathrm{GL} }\nolimits}
\newcommand{\Sp}{\mathop{\mathrm{Sp} }\nolimits}
\newcommand{\GSp}{\mathop{\mathrm{GSp} }\nolimits}
\newcommand{\PGL}{\mathop{\mathrm{PGL} }\nolimits}

\newcommand\F{\mathbb{F}}

\newcommand{\Aut}{\mathop{\mathrm{Aut} }\nolimits}

\newcommand{\WD}{\mathop{\mathrm{WD} }\nolimits}

\newcommand{\rec}{\mathop{\mathrm{rec} }\nolimits}
\newcommand\SL{\mathrm{SL}}
\newcommand\SO{\mathrm{SO}}
\newcommand\PSp{\mathrm{PSp}}
\newcommand\PSL{\mathrm{PSL}}
\newcommand\PGSp{\mathrm{PGSp}}

\newcommand\End{\mathrm{End}}

\begin{document}

\title{Automorphic Galois representations and the inverse Galois problem}

\author{Sara Arias-de-Reyna\footnote{Universit\'{e} du Luxembourg,
Facult\'{e} des Sciences, de la Technologie et de la Communication,
6, rue Richard Coudenhove-Kalergi, L-1359 Luxembourg, Luxembourg,
sara.ariasdereyna@uni.lu}}

\date{}

\maketitle

\begin{abstract} A strategy to address the inverse Galois problem over $\Q$ consists of exploiting the knowledge of  Galois representations attached to certain automorphic forms. More precisely, if such forms are carefully chosen, they provide compatible systems of Galois representations satisfying some desired properties, e.g. properties that reflect on the image of the members of the system.
In this article we survey some results obtained using this strategy.

MSC (2010): 11F80 (Galois representations);
12F12 (Inverse Galois theory).

\end{abstract}

\section{Introduction}

The motivation for the subject of this survey comes from Galois theory. 
Let $L/K$ be a field extension which is normal and separable. To this extension one can attach a group, namely the group of field automorphisms of $L$ fixing $K$, which is denoted as $\Gal(L/K)$.
The main result of Galois theory, which is usually covered in the program of any Bachelor's degree in Mathematics, can be stated as follows:

\begin{theorem}[Galois]
Let $L/K$ be a finite, normal, separable field extension. Then there is the following bijective correspondence between the sets:

\begin{equation*}
\begin{array}{ccc}
\left\{\begin{array}{c} E \text{ field} \\ K\subseteq E\subseteq L\end{array}\right\} & \longleftrightarrow &
                     \left\{\begin{array}{c} H\subseteq\Gal(L/K) \\ \text{subgroup}\end{array}\right\},\\
[2.3em]
               E & \longmapsto & \mathrm{Gal}(L/E)\\
                 L^H & \longmapsfrom& L\end{array}
                     \end{equation*}

\end{theorem}

Usually the students are asked exercises of the following type: Given some finite field extension $L/\Q$ which is normal, compute the Galois group $\mathrm{Gal}(L/\Q)$ attached to it. But, one may also ask the inverse question (hence the name inverse Galois problem): Given a finite group $G$, find a finite, normal extension $L/\Q$ with $\Gal(L/\Q)\simeq G$. This is not a question one usually expects a student to solve! In fact, there are (many) groups $G$ for which it is not even known if such a field extension exists.

\begin{prob}[Inverse Galois Problem] Let $G$ be a finite group. Does there exist a Galois extension $L/\Q$ such that $\Gal(L/\Q)\simeq G$?
\end{prob}

The first mathematician that addressed this problem was D.~Hilbert. In his paper \cite{Hilbert} he proves his famous Irreducibility Theorem, and applies it to show that, for all $n\in \mathbb{N}$, the symmetric group $S_n$ and the alternating group $A_n$ occur as Galois groups over the rationals. Since then, many mathematicians have thought about the inverse Galois problem, and in fact it is now solved (affirmatively) for many (families of) finite groups $G$. For instance, let us mention the result of Shafarevich that all solvable groups occur as Galois groups over the rationals (see \cite{Cohomology_of_number_fields} for a detailed explanation of the proof). 
However, it is still not known if the answer is affirmative for every finite group $G$, and as far as I know, there is no general strategy that addresses all finite groups at once.  An account of the different techniques used to address the problem can be found in \cite{Topics}.
%The reader can find many details in \cite{MM}.

Let $K$ be a field, and let us fix a separable closure $K_{\rm sep}$. There is a way to group together all the Galois groups of finite Galois extensions $L/K$ contained in $K_{\rm sep}$, namely the \emph{absolute Galois group} of $K$. It is defined as the inverse limit
\begin{equation*}
 G_{K}:=\Gal(K_{\rm sep}/K)=\lim_{\longleftarrow \atop {L/K\atop \text{finite Galois }}}\Gal(L/K).
\end{equation*}
This group is a profinite group, and as such is endowed with a topology, called the \emph{Krull topology}, which makes it a Hausdorff, compact and totally disconnected group. A very natural question to ask is what information on the field $K$ is encoded in the topological group $G_K$. In this connection, a celebrated result of Neukirch, Iwasawa, Uchida and Ikeda establishes that, if $K_1, K_2$ are two finite extensions of $\Q$ contained in a fixed algebraic closure $\overline{\Q}$ such that $G_{K_1}\simeq G_{K_2}$, then $K_1$ and $K_2$ are conjugated by some element in $G_{\Q}$ (cf.~\cite{Uchida1976}, \cite{Ikeda}). Let us note, however, that we cannot replace $\Q$ by any field. For example, the analogous statement does not hold when the base field is $\Q_p$, cf.~\cite{Yamagata} and \cite{JardenRitter}.
Thus, we see that the absolute Galois group of $\Q$ encodes a wealth of information about the arithmetic of number fields. In this context, the inverse Galois problem can be formulated as the question of determining which finite groups occur as quotient groups of $G_{\Q}$.

A natural way to study $G_{\Q}$ is to consider its representations, that it, the continuous group morphisms $G_{\Q}\rightarrow \GL_m(k)$, where $k$ is a topological field and $m\in \mathbb{N}$. Such a representation will be called a \emph{Galois representation}. Let us assume that $k$ is a finite field, endowed with the discrete topology, and let
\begin{equation*}
 \rho:G_{\Q}\rightarrow \GL_m(k)
\end{equation*} be a Galois representation. Since the set $\{\mathrm{Id}\}$ is open in $\GL_m(k)$, we obtain that $\ker\rho\subset G_{\Q}$ is an open subgroup. In other words, there exists a finite Galois extension $K/\Q$ such that $\ker\rho= G_K$. Therefore\
\begin{equation*}
 \mathrm{Im}\rho\simeq G_{\Q}/\ker\rho\simeq G_{\Q}/G_{K}\simeq \Gal(K/\Q).
\end{equation*}

This reasoning shows that, whenever we are given a Galois representation of $G_{\Q}$ over a finite field $k$, we obtain a realisation of $\mathrm{Im}\rho\subset \GL_m(k)$ as a Galois group over $\Q$. In this way, any source of Galois representations provides us with a strategy to address the inverse Galois problem for the subgroups of $\GL_m(k)$ that occur as images thereof.

Geometry provides us with many objects endowed with an action of the absolute Galois group of the rationals, thus giving rise to such Galois representations. One classical example is the group of $\overline{\Q}$-defined $\ell$-torsion points of an elliptic curve $E$ defined over $\Q$. We will treat this example in Section \ref{sec:Classical}. In this survey we will be interested in (compatible systems of) Galois representations arising from automorphic representations. In Section \ref{sec:Automorphic} we will describe Galois representations attached to an automorphic representation $\pi$ which satisfies several technical conditions. The statements of the most recent results (to the best of my knowledge) on the inverse Galois problem obtained by means of compatible systems of Galois representations attached to automorphic representations can be found in Section \ref{sec:Ingredients}, together with some ideas about their proofs.

A remarkable feature of this method is that, in addition, one obtains some control of the ramification of the Galois extension that is produced. Namely, it will only be ramified at the residual characteristic and at a finite set of auxiliary primes, that usually one is allowed to choose (inside some positive density set of primes). This will be highlighted in the statements below. 

\smallskip

\textbf{Acknowledgements:} This article is an expanded version of the plenary lecture I delivered at the conference \emph{Quintas Jornadas de Teor\'ia de N\'umeros} (July 2013). I would like to thank the scientific committee for giving me the oportunity to participate in this conference, and the organising committee for their excellent work. I also want to thank Gabor Wiese for his remarks and suggestions on a previous version of this article.

\section{Some classical cases}\label{sec:Classical}

In this section we revisit some classical examples of Galois representations attached to geometric objects. We begin with the Galois representations attached to the torsion points of elliptic curves, and later we will see them as a particular case of Galois representations attached to modular forms.

\subsection{Elliptic curves}

An elliptic curve is a genus one curve, endowed with a distinguished base point. Every elliptic curve $E$ can be described by means of a Weierstrass equations, that is,  an affine equation of the form
\begin{equation*}
 y^2 + a_1xy + a_3 y=x^3 + a_2 x^2 + a_4 x + a_6
\end{equation*}
where the coefficients $a_1, \dots, a_6$ lie in some field $K$. The most significant property of elliptic curves is that the set of points of $E$ (defined over some field extension $L/K$) can be endowed with a commutative group structure, where the neutral element is the distinguished base point.

Let $E/\mathbb{Q}$ be an elliptic curve and let $\ell$ be a prime number. We can consider the subgroup $ E[\ell](\overline{\Q})$ of $E(\overline{\Q})$ consisting of $\ell$-torsion points. This group is isomorphic to the product of two copies of $\mathbb{F}_{\ell}$. Moreover, since the elliptic curve is defined over $\mathbb{Q}$, the absolute Galois group $G_{\Q}$ acts naturally on the set of $\overline{\Q}$-defined points of $E$, and this action restricts to $ E[\ell](\overline{\Q})$. We obtain thus a Galois representation
\begin{equation*}
\overline{\rho}_{E, \ell}:G_{\Q}\rightarrow \Aut(E[\ell](\overline{\Q}))\simeq \GL_2(\mathbb{F}_{\ell}).
\end{equation*}

As explained in the introduction, the image of $\overline{\rho}_{E, \ell}$ can be realised as a Galois group over $\Q$. This brings forward the question of determining the image of such a Galois representation. In this context, there is a classical result by J.~P.~Serre from the seventies (\cite{Proprietes}, Th\'eor\`eme 2).

\begin{theorem}[Serre] Let $E/\Q$ be an elliptic curve without complex multiplication over $\overline{\Q}$. Then the representation $\overline{\rho}_{E, \ell}$ is surjective for all except finitely many primes $\ell$.
\end{theorem}

We can immediately conclude that $\GL_2(\mathbb{F}_{\ell})$ can be realised as a Galois group over $\Q$ for all except finitely many primes $\ell$. However, we can do even better by picking a particular elliptic curve and analysing the Galois representations attached to it.

\begin{example}\label{ex:37}
Let $E/\Q$ be the elliptic curve defined by the Weierstrass equation
\begin{equation*}
 y^2 + y=x^3-x.
\end{equation*}

This curve is labelled 37A in \cite{CremonaTables}, and it has the property that $\overline{\rho}_{E, \ell}$ is surjective \emph{for all primes $\ell$} (see \cite{Proprietes}, Example 5.5.6). Therefore we obtain that $\GL_2(\mathbb{F}_{\ell})$ occurs as the Galois group of a finite Galois extension $K/\Q$. Moreover, we have additional information on the ramification of $K/\Q$; namely, it ramifies only at $37$ (which is the conductor of $E$) and $\ell$.
\end{example}

The next situation we want to analyse is that of Galois representations attached to modular forms. Let us recall that
modular forms are holomorphic functions defined on the complex upper half plane, which satisfy certain symmetry relations. We will not recall here the details of the definition (see e.g.~\cite{Diamond-Shurman} for a complete treatment focusing on the relationship with arithmetic geometry). These objects, of complex-analytic nature, play a central role in number theory. At the core of this relationship is the fact that one can attach Galois representations of $G_{\Q}$ to them.
More precisely, let $f$ be a cuspidal modular form of weight $k\geq 2$, conductor $N$ and character $\psi$ (in short: $f\in S_k(N, \psi)$), which is a normalised Hecke eigenform. We may write the Fourier expansion of $f$ as $f(z)=\sum_{n\geq 1}a_n q^n$, where $q=e^{2\pi i z}$.
A first remark is that the \emph{coefficient field} $\Q_f=\Q(\{a_n:\gcd(n, N)=1\})$ is a number field. Denote by $\mathcal{O}_{\Q_f}$ its ring of integers. By a result of
Deligne (cf.~\cite{De71}), for each prime $\lambda$ of $\mathcal{O}_{\Q_f}$ there exists a (continuous) Galois representation
 \begin{equation*}
  \rho_{f, \lambda}: G_{\Q}\rightarrow \GL_2(\mathcal{O}_{\overline{\Q}_{f,\lambda}}),
 \end{equation*}
related to $f$, where $\Q_{f,\lambda}$ denotes the completion of $\Q_f$ at the prime $\lambda$, $\overline{\Q}_{f,\lambda}$ an algebraic closure thereof and $\mathcal{O}_{\overline{\Q}_{f,\lambda}}$ is the valuation ring of $\overline{\Q}_{f,\lambda}$. Here the topology considered on $\GL_2(\mathcal{O}_{\overline{\Q}_{f,\lambda}})$ is the one induced by the $\ell$-adic valuation.

The relationship between $\rho_{f, \lambda}$ and $f$ is the following. First, $\rho_{f, \lambda}$ is unramified outside $N\ell$. Moreover, for each $p\nmid N\ell$, we can consider the image under $\rho_{f, \lambda}$ of a lift $\mathrm{Frob}_p$ of a Frobenius element at $p$ (this is well defined because $\rho_{f, \lambda}$ is unramified at $p$). Then the characteristic polynomial of $\rho_{f, \lambda}(\mathrm{Frob}_p)$ equals $T^2 - a_p T + \psi(p) p^{k-1}$.

We may compose each $\rho_{f, \lambda}$ with the reduction modulo the maximal ideal of $\mathcal{O}_{\overline{\Q}_{f,\lambda}}$, and we obtain a (residual) representation
 \begin{equation*}
  \overline{\rho}_{f, \lambda}: G_{\Q}\rightarrow \GL_2(\kappa(\overline{\Q}_{f,\lambda}))\simeq  \GL_2(\overline{\F}_{\ell}),
 \end{equation*}
where $\ell$ is the rational prime below $\lambda$.
One of the main recent achievements in number theory has been the proof of Serre's Modularity Conjecture, which says that every Galois representation $ \overline{\rho}_{\ell}: G_{\Q}\rightarrow \GL_2(\overline{\F}_{\ell})$ which is odd and irreducible is actually isomorphic to $\overline{\rho}_{f, \lambda}$ for some modular form $f$ and some prime $\lambda$ as above.

In this survey we are interested in the image of $\overline{\rho}_{f, \lambda}$. These images have been studied by K.~Ribet (cf.~\cite{Ribet75}, \cite{Ribet85}). One first remark is that, when $\rho_{f, \lambda}$ is absolutely irreducible, then it can be conjugated (inside $\GL_2(\mathcal{O}_{\overline{\Q}_{f,\lambda}})$) so that its image is contained in $\GL_2(\mathcal{O}_{\Q_{f,\lambda}})$. Therefore, in this case we can assume that $\overline{\rho}_{f, \lambda}: G_{\Q}\rightarrow \GL_2(\kappa(\Q_{f,\lambda}))$, where $\kappa(\Q_{f,\lambda})$ denotes the residue field of $\Q_{f,\lambda}$.

To state Ribet's result, we first introduce two more number fields related to $f$. The first one is the \emph{twist invariant coefficient field of $f$}, which is the subfield of the coefficient field of $f$ defined as $F_f:=\Q(\{a_n^2/\psi(n): \gcd(n, N)=1\})$. The second field, which is a finite abelian extension of $\Q$,  is the subfield $K_f$ of  $\overline{\Q}$  fixed by all inner twists of $f$ (see \cite{DiWi} for details).

\begin{theorem}[Ribet]\label{theorem:Ribet} Let $f=\sum_{n\geq 1}a_n q^n\in S_k(N, \psi)$ be a normalised cuspidal Hecke eigenform. Assume $f$ does not have complex multiplication. Then for all except finitely many prime ideals $\lambda$ of $\Q_{f}$,
\begin{equation*}
 \overline{\rho}_{f, \lambda}(G_{K_f})=\left\{g\in \GL_2\left(\kappa\left(F_{f, \lambda'}\right)\right):\det(g)\in (\mathbb{F}_{\ell}^{\times})^{k-1}\right\},
\end{equation*}
where $\lambda'$ is the ideal of $\mathcal{O}_{F_f}$ below $\lambda$ and $\ell$ is the rational prime below $\lambda$.
\end{theorem}

This result suggests that we look at the representation $\overline{\rho}^{\mathrm{proj}}_{f, \lambda}$ obtained by composing $\overline{\rho}_{f, \lambda}$ with the projection  map $\GL_2(\kappa(\Q_{f,\lambda}))\rightarrow \PGL_2(\kappa(\Q_{f, \lambda}))$.

More precisely, let $k, r$ be integers greater than or equal to $1$. Consider the set 
\begin{equation*} \mathcal{A}:=\{A\in \GL_2(\mathbb{F}_{\ell^r}): \det A\in (\mathbb{F}_{\ell}^{\times})^{k-1}\},\end{equation*} and let $\mathcal{A}^{\mathrm{proj}}$ be its projection under the map $\GL_2(\mathbb{F}_{\ell^r})\rightarrow \PGL_2(\mathbb{F}_{\ell^r})$. Then if $k$ is odd, we have $\mathcal{A}^{\mathrm{proj}}=\PSL_2(\mathbb{F}_{\ell^r})$, and if $k$ is even, we have $\mathcal{A}^{\mathrm{proj}}=\PGL_2(\mathbb{F}_{\ell^r})$ whenever $r$ is odd and $\mathcal{A}^{\mathrm{proj}}=\PSL_2(\mathbb{F}_{\ell^r})$ whenever $r$ is even.

In any case it follows that, for $f$ as above, the image of $\overline{\rho}^{\mathrm{proj}}_{f, \lambda}$ equals $\PSL_2\left(\kappa\left(F_{f, \lambda'}\right)\right)$ or $\PGL_2\left(\kappa\left(F_{f, \lambda'}\right)\right)$ for all except finitely many primes $\lambda$ of $\mathcal{O}_{\Q_f}$.

A remarkable difference with the situation arising from elliptic curves is that we obtain realisations of linear groups over fields whose cardinality is not necessarily a prime number. In Example \ref{ex:37}, we used an elliptic curve to obtain realisations of the members of the family $\{\GL_2(\mathbb{\F}_{\ell})\}_{\ell}$. However, now we have two parameters, namely the prime $\ell$ and the exponent $r$. If we pick a modular form as above, we will obtain realisations of members of one of the families $\{\PSL_2(\mathbb{F}_{\ell^r})\}_{\ell, r}$ or  $\{\PGL_2(\mathbb{F}_{\ell^r})\}_{\ell, r}$, and the parameter $r$ \emph{depends on $f$ and $\ell$}.

\begin{example}[Ribet, 1975]\label{ex:Ribet}
Let $f\in S_{24}(1)$ be a normalised Hecke eigenform of level $1$. The field of coefficients $\Q_f=\Q(\sqrt{144169})$ equals $F_f$; so we can expect to obtain realisations of $\PSL_2(\mathbb{F}_{\ell^2})$ when $\ell$ is inert in $\Q_f$ and $\PGL_2(\mathbb{F}_{\ell})$ when $\ell$ splits in $\Q_f$. Indeed, let $\ell$ be a prime different from $2$, $3$ and $47$. Then $f$ provides a realisation of $\PGL_2(\mathbb{F}_{\ell})$ if $144169$ is a square modulo $\ell$ and a realisation of $\PSL_2(\mathbb{F}_{\ell^2})$ if $144169$ is not a square modulo $\ell$. Moreover, the corresponding Galois extension $K/\Q$ with desired Galois group is unramified outside $\ell$.
\end{example}

\section{Compatible systems and the inverse Galois
problem}\label{sec:CompatibleSystems}

The examples of the previous section suggest that, instead of considering isolated Galois representations $\overline{\rho}:G_{\Q}\rightarrow \GL_n(\overline{\mathbb{F}}_{\ell})$ for a fixed prime $\ell$, it is a good idea to look at a system of Galois representations $(\overline{\rho}_{\ell})_{\ell}$, where $\ell$ runs through the prime numbers. The notion of \emph{(strictly) compatible system of Galois representations} already appears in \cite{SerreAbelian}. We recall the definition below.

\begin{definition}\label{defi:CompatibleSystem} Let $n\in \mathbb{N}$ and let $F$ be a number field. A {\em compatible system $\rho_\bullet = (\rho_\lambda)_\lambda$
of $m$-dimensional representations of $G_{F}$} consists of the following data:
\begin{itemize}
\item A number field $L$.
\item A finite set $S$ of primes of $F$.
\item For each prime $\mathfrak{p}\not\in S$, a monic polynomial $P_\mathfrak{p}(X) \in \mathcal{O}_L[X]$
(with $\mathcal{O}_L$ the ring of integers of~$L$).
\item For each finite place $\lambda$ of~$L$
(together with fixed embeddings $L \hookrightarrow L_\lambda \hookrightarrow \overline{L}_\lambda$) a continuous Galois representation
$$ \rho_\lambda: G_F \to \GL_m(\overline{L}_\lambda)$$
such that $\rho_\lambda$ is unramified outside $S \cup S_{\ell}$
(where $\ell$ is the rational prime below~$\lambda$ and $S_{\ell}$ is the set of primes of $F$ above $\ell$)
and such that for all~$\mathfrak{p} \not\in S \cup S_{\ell}$ the characteristic
polynomial of $\rho_\lambda(\mathrm{Frob}_{\mathfrak{p}})$ is equal to $P_{\mathfrak{p}}(X)$.
\end{itemize}
\end{definition}

In our context, the main question to ask about a compatible system is the following: If we know that $\rho_{\lambda}$ satisfies some property (A), does it follow that $\rho_{\lambda'}$ also satisfies (A) for another prime $\lambda'$ of $L$? In other words, what properties ``propagate'' through a compatible system? The idea that the property of ``being attached to a modular form'' propagates through such a system lies at the core of the proof of the Taniyama-Shimura conjecture by A.~Wiles and R.~Taylor (which implies Fermat's Last Theorem), and also of the proof of Serre's Modularity Conjecture.

In this section we are interested in the relationship between the images of the members $\rho_{\lambda}$ of a compatible system. An example of such a relationship is the following: if $\rho_{\lambda}$, $\rho_{\lambda'}$ are two semisimple representations belonging to a compatible system, then the image of $\rho_{\lambda}$ is abelian if and only if the image of $\rho_{\lambda'}$ is abelian (see \cite{SerreAbelian} and \cite{Henniart1982}).

The case of compatible systems of Galois representations attached to the Tate module of abelian varieties has received particular attention. Let $A/F$ be an $n$-dimensional abelian variety, and assume that
\begin{equation*}(\rho_{A, \ell}:G_{F}\rightarrow \mathrm{GL}(V_{\mathbb{Q}_{\ell}})\simeq \mathrm{GL}_{2n}(\mathbb{Q}_{\ell}))_{\ell}\end{equation*}
is the compatible system of Galois representations attached to the $\ell$-adic Tate module $T_{\ell}$ of $A$ (where as usual
$V_{\mathbb{Q}_{\ell}}=\mathbb{Q}_{\ell}\otimes_{\mathbb{Z}_{\ell}}T_{\ell})$). To what extent does the image of $\rho_{\ell}$ depend on $\ell$?
There are several ways to phrase this question in a precise way. For example, define the \emph{algebraic monodromy group} at $\ell$, $G_{\ell}$, as the Zarisky closure of $\rho_{\ell}(G_F)$ inside the algebraic group $\GL_{2n, \mathbb{Q}_{\ell}}$, and let $G_{\ell}^0$ be the connected component of $G_{\ell}$. In this connection, the Mumford-Tate conjecture predicts the existence of an algebraic group $G\subset \GL_{2n, \mathbb{Q}}$ such that, for all $\ell$, $G_{\ell}^0\simeq \mathbb{Q}_{\ell}\times_{\mathbb{Q}}G$ (see \cite{Serre:Iyanaga}, Conjecture C.3.3 for a precise formulation). By work of J.~P.~Serre it is known that the (finite) group of connected components $G_{\ell}/G_{\ell}^0$ is independent of $\ell$ (see \cite{Serre:Course}, 2.2.3).

There are many partial results in this direction. In particular cases, the conjecture is known to hold (for example, when $\dim A=1$ cf.~\cite{SerreAbelian} and \cite{Proprietes}. For higher dimension, when $\End_{\overline{\Q}}(A)=\mathbb{Z}$ and $n=2$ or odd the conjecture holds with $G=\GSp_{2n, \Q}$; cf.~\cite{Serre:Course}, 2.2.8). In the general case, Serre has proved that the rank of $G_{\ell}$ is independent of $\ell$ \cite{Serre:RibetLetter}. More partial results can be found in \cite{LarsenPink1992}, \cite{Hui2013}.

Another question is how close $\rho_{\ell}(G_{F})$ is to its Zarisky closure  $G_{\ell}$ in $\GL_{2g, \Q_{\ell}}$. For results in this direction the reader is referred to \cite{Larsen1995} and \cite{HuiLarsen}.

A particular case, which is of interest to us (cf. Section \ref{sec:Ingredients}), is proved by C.~Hall in \cite{Hall2011}. Let $A/F$ be an $n$-dimensional abelian variety which is principally polarised and with $\mathrm{End}_{\overline{\Q}}(A)=\mathbb{Z}$. Assume that there exists a prime $\mathfrak{p}$ of $F$ such that the reduction of $A$ at $\mathfrak{p}$ is semistable of toric dimension $1$. Then
there exists a constant $M$ such that, for all primes $\ell\geq M$, the image of the mod $\ell$ Galois representation $\overline{\rho}_{A, \ell}:G_F\rightarrow \GL_{2n}(\mathbb{F}_{\ell})$ coincides with $\mathrm{GSp}_{2n}(\mathbb{F}_{\ell})$. As a consequence, it follows that
$A$ satisfies the Mumford-Tate conjecture; more precisely, the corresponding algebraic group is $\mathrm{GSp}_{2n, \mathbb{Q}}$. The proof of this result relies heavily on the fact that the existence of the prime $\mathfrak{p}$ implies that the image under $\overline{\rho}_{A, \ell}$ of the inertia group at $\mathfrak{p}$ contains a transvection.

For the applications to the inverse Galois problem, we will be interested in Galois representations taking values in linear groups over finite fields. For the rest of the section, we focus on symplectic groups $\GSp_{2n}$ for simplicity. Note that $\GSp_2=\GL_2$ and $\Sp_2=\SL_2$, so in the case of dimension $1$ we are in the situation explained in Section \ref{sec:Classical}. Consider the following setup:

\begin{setup}\label{setup:symplectic} Let $\rho_{\bullet}=(\rho_{\lambda})_{\lambda}$ be a $2n$-dimensional compatible system of Galois representations of $G_{\mathbb{Q}}$ as in Definition \ref{defi:CompatibleSystem}, such that for all $\lambda$, $\rho_{\lambda}:G_{\Q}\rightarrow \GSp_{2n}(\overline{L}_{\lambda})$  for some number field $L$ (we will say that such a system is \emph{symplectic}).
\end{setup}

Note that each of the $\rho_{\lambda}$ is defined over a finite extension of $L_{\lambda}$ inside $\overline{L}_{\lambda}$. Moreover, we can conjugate each $\rho_{\lambda}$ to take values inside the ring of integers of this finite extension of $L_{\lambda}$, and further reduce it modulo $\lambda$, obtaining a residual representation $\overline{\rho}_{\lambda}$. When $\overline{\rho}_{\lambda}$ is absolutely irreducible, then $\rho_{\lambda}$ can be defined over $L_{\lambda}$, and therefore $\overline{\rho}_{\lambda}$ takes values inside $\mathrm{GSp}_{2n}(\kappa(L_{\lambda}))$, where $\kappa(L_{\lambda})$ denotes the residue field of $L_{\lambda}$. Recall the motivating example in Section \ref{sec:Classical} of compatible systems attached to modular forms. In this example, the field $L$ can be taken to be the coefficient field $\mathbb{Q}_f$. Like in the case of compatible systems attached to modular forms, it will be convenient to consider the composition $\overline{\rho}_{\lambda}^{\mathrm{proj}}$ of $\overline{\rho}_{\lambda}$ with the natural projection $\mathrm{GSp}_{2n}(\kappa(L_{\lambda}))\rightarrow \mathrm{PGSp}_{2n}(\kappa(L_{\lambda}))$. In what follows, we focus on obtaining realisations of groups in one of the families $\{\mathrm{PSp}_{2n}(\mathbb{F}_{\ell^r})\}_{\ell, r}$ or $\{\mathrm{PGSp}_{2n}(\mathbb{F}_{\ell^r})\}_{\ell, r}$.

Assume that we are given a compatible system of Galois representations as in Set-up \ref{setup:symplectic} such that all $\rho_{\lambda}$ are residually absolutely irreducible. We obtain a system \begin{equation*}(\overline{\rho}_{\lambda}:G_{\mathbb{Q}}\rightarrow \GSp_{2n}(\kappa(L_{\lambda})))_{\lambda}.\end{equation*} For each prime $\lambda$ of $L$, $\kappa(L_{\lambda})\simeq \mathbb{F}_{\ell^{r(\lambda)}}$ for some integer $r(\lambda)$, which actually changes with $\lambda$! If we want to realise the family of groups $\{\PSp_{2n}(\mathbb{F}_{\ell^r})\}_{\ell}$ for a fixed exponent $r$, it is clear that one compatible system will not suffice for our purposes (unless we are interested in $r=1$ and we have $L=\mathbb{Q}$). This phenomenon already appeared in Section \ref{sec:Classical} in the case of compatible systems attached to modular forms.

The strategy to obtain Galois realisations will proceed as follows. We want to construct a compatible system of Galois representations $\rho_{\bullet}$ as in Set-up \ref{setup:symplectic}, such that the $\rho_{\lambda}$ are absolutely irreducible, and such that the images of the corresponding representations $\overline{\rho}_{\lambda}$ are large in some sense which does not depend on $\lambda$. More precisely, we will say that the image of a representation $\overline{\rho}_{\lambda}:G_{\Q}\rightarrow \GSp_{2n}(\kappa(L_{\lambda}))$ is \emph{huge} if it contains a conjugate (inside $\GSp_{2n}(\overline{\mathbb{F}}_{\ell})$) of $\mathrm{Sp}_{2n}(\mathbb{F}_{\ell})$ (where $\ell$ is the prime below $\lambda$). A group theoretical reasoning shows that if $\overline{\rho}_{\lambda}$ has huge image, then the image of $\overline{\rho}^{\mathrm{proj}}_{\lambda}$ equals $\PGSp_{2n}(\mathbb{F}_{\ell^r})$ or $\PSp_{2n}(\mathbb{F}_{\ell^r})$ for some integer $r$ (cf.~Corollary 5.7 of \cite{ArDiWi1}). Moreover, we will have to find some conditions to control the exponent $r$.

\smallskip

The presence of these two parameters, $\ell$ and $r$, gives rise to two different approaches to obtain results on the inverse Galois problem:

\begin{itemize}
 \item \textbf{Vertical Direction:} Fix a prime number $\ell$. Obtain realisations of $\mathrm{PSp}_{2n}(\mathbb{F}_{\ell^r})$ (resp. $\mathrm{PGSp}_{2n}(\mathbb{F}_{\ell^r})$)  for all $r\in \mathbb{N}$.

 \item \textbf{Horizontal Direction:} Fix a natural number $r\geq 1$. Obtain realisations of $\mathrm{PSp}_{2n}(\mathbb{F}_{\ell^r})$ (resp. $\mathrm{PGSp}_{2n}(\mathbb{F}_{\ell^r})$)  for all primes $\ell$.
\end{itemize}

This nomenclature stems from the following representation: Place in a graphic the groups in the family $\mathrm{PSp}_{2n}(\mathbb{F}_{\ell^r})$ (resp. $\mathrm{PGSp}_{2n}(\mathbb{F}_{\ell^r})$) that are realised as Galois groups over $\Q$ by displaying in the $x$-axis the prime $\ell$ and in the $y$-axis the exponent $r$, and drawing a dot whenever the group $\mathrm{PSp}_{2n}(\mathbb{F}_{\ell^r})$ (resp. $\mathrm{PGSp}_{2n}(\mathbb{F}_{\ell^r})$) is realised as a Galois group over $\Q$ (see \cite{DiWi} for such a graphic when $n=1$).

By exploiting the compatible systems of Galois representations attached to modular forms, the following results have been proved in the vertical direction (see Theorem 1.1 of \cite{Wiese2008}) and in the horizontal direction (see Theorem 1.1 of \cite{DiWi}).

\begin{theorem}[Wiese]\label{theorem:Wiese} Let $\ell$ be a prime number. There exist infinitely many natural numbers $r$ such that $\PSL_2(\mathbb{F}_{\ell^r})$ occurs as the Galois group of a finite Galois extension $K/\mathbb{Q}$, which is unramified outside $\ell$ and an auxiliary prime $q$. \end{theorem}

\begin{theorem}[Dieulefait, Wiese]\label{theorem:DieulefaitWiese} Let $r\in \mathbb{N}$.

\begin{enumerate}
  \item There exists a positive density set of primes $\ell$ such that $\mathrm{PSL}_2(\mathbb{F}_{\ell^r})$ occurs as the Galois group of a finite Galois extension $K/\mathbb{Q}$, which is unramified outside $\ell$ and two (resp. three) auxiliary primes if $n$ is even (resp. odd).
  \item Assume that $r$ is odd. There exists a positive density set of primes $\ell$ such that $\mathrm{PGL}_2(\mathbb{F}_{\ell^r})$ occurs as the Galois group of a finite Galois extension $K/\mathbb{Q}$, which is unramified outside $\ell$ and two auxiliary primes.
 \end{enumerate}
\end{theorem}

Let us look more closely at the approach in the horizontal direction. We fix a natural number $r$, and we want to realise $\mathrm{PSL}_2(\mathbb{F}_{\ell^r})$ or $\mathrm{PGL}_{2}(\mathbb{F}_{\ell^r})$ as Galois groups over $\Q$ for as many primes $\ell$ as we can. From the remarks above, it is clear that a single modular form will not suffice to realise $\mathrm{PSL}_2(\mathbb{F}_{\ell^r})$ for \emph{all} $\ell$. However, nothing prevents us from looking at several modular forms. In fact, Serre's Modularity Conjecture, which is now a theorem, tells us that every irreducible, odd Galois representation $\overline{\rho}:G_{\Q}\rightarrow \mathrm{GL}_2(\overline{\mathbb{F}}_{\ell})$ is attached to some modular form $f$. As a consequence, any realisation of $\mathrm{PSL}_2(\mathbb{F}_{\ell^r})$ as the Galois group of a finite Galois extension $K/\mathbb{Q}$ with $K$  imaginary can be obtained through this method (cf.~Proposition 1.2 of \cite{DiWi}). By making use of Theorem \ref{theorem:Ribet}, we know that for a normalised Hecke eigenform without complex multiplication, the image of $\overline{\rho}_{f, \lambda}$ is huge for all except finitely many prime ideals $\lambda$ of $\mathbb{Q}_f$, and thus the image of $\overline{\rho}^{\mathrm{proj}}_{f, \lambda}$ is isomorphic to $\mathrm{PSL}_2(\mathbb{F}_{\ell^r})$ or $\mathrm{PGL}_2(\mathbb{F}_{\ell^r})$. The main obstacle here is to obtain some control on the exponent $r$. Under additional conditions, the field $\mathbb{F}_{\ell^r}$ coincides with $\kappa(\Q_{f, \lambda})$, reducing the problem to the analysis of $\mathbb{Q}_f$. But this is not a minor issue! Very little is known about these fields (although one can always compute them for any given modular form $f$). When the level of $f$ is $1$, there is a strong conjecture in this connection, namely Maeda's conjecture, stating that the degree $d_f=[\mathbb{Q}_f:\mathbb{Q}]$ should equal the dimension of $S_k(1)$ as a complex vector space ($k$ being the \emph{weight} of $f$) and the Galois group of the normal closure of $\mathbb{Q}_f/\mathbb{Q}$ is equal to the symmetric group $S_{d_f}$. Assuming this conjecture, one can improve Theorem \ref{theorem:DieulefaitWiese} as follows (cf.~Theorem 1.1 of \cite{Wiese2013}).

\begin{theorem}[Wiese] Assume Maeda's Conjecture holds. Let $r\in \mathbb{N}$.
Assume that $r$ is even (resp.~odd). There exists a density 1 set of primes $\ell$ such that $\mathrm{PSL}_2(\mathbb{F}_{\ell^r})$ (resp.~$\mathrm{PGL}_2(\mathbb{F}_{\ell^r})$)  occurs as the Galois group of a finite Galois extension $K/\mathbb{Q}$, which is unramified outside $\ell$.
\end{theorem}

\section{Galois representations attached to automorphic forms}\label{sec:Automorphic}

In order to use the strategy outlined in the
previous section to obtain results on the inverse Galois problem, we
first need to find a source of compatible systems of Galois
representations of $G_{\mathbb{Q}}$. As discussed in Section
\ref{sec:Classical}, elliptic curves defined over $\mathbb{Q}$ (and,
analogously, abelian varieties of higher dimension which are defined
over $\mathbb{Q}$) provide us with such systems, and, more
generally, classical modular forms give rise to such systems. Both
of these examples can be encompassed in the general framework
provided by the Langlands conjectures. More precisely, given an
automorphic representation $\pi$ (which is \emph{algebraic} in some
precise sense) for an arbitrary connected reductive group $G$ over $\Q$, one hopes that there exists a compatible system
of Galois representations $(\rho_{\bullet}(\pi))$ attached to it, where $\rho_{\lambda}(\pi)$
takes values in the $\overline{\Q}_{\ell}$-points of a certain
algebraic group attached to $G$ (namely the Langlands dual of $G$). Conjecturally, then, we have many
compatible systems of Galois representations, which builds up the
hope of eventually applying the strategy described in the previous
section to realise many linear groups as Galois groups over the rationals.

There are several cases when these conjectures are known to hold. Recently, there has been  a breakthrough in this
connection due to P.~Scholze \cite{Scholze} and M.~Harris, K.-W.~Lan, R.~Taylor,
J.~Thorne. Namely, they attach compatible systems of Galois
representations to regular, L-algebraic cuspidal automorphic
representations of $\mathrm{GL}_m(\mathbb{A}_F)$, where $F$ is a
totally real or a CM number field. 

However, in this section we will recall a less recent result, due to L.~Clozel, R.~Kottwitz, M.~Harris, R.~Taylor and several
others, which is more restrictive, since it deals with RAESDC
(regular, algebraic, essentially self-dual, cuspidal) automorphic
representations.  We will not recall here all definitions (the reader can look
them up in \cite{BLGGT}), but we will try to give some explanations. 

Let $\mathbb{A}_{\Q}$ denote the ring of adeles of $\Q$. We consider so-called irreducible admissible representations $\pi$ of $\GL_{m}(\mathbb{A}_{\Q})$. In fact, $\pi$ is not literally a representation of the group $\GL_{m}(\mathbb{A}_{\Q})$ into some vector space. The interested reader can look at the details in \cite{Bump}. In this survey, we will treat them as black boxes, focusing rather on the compatible systems of Galois representations that they give rise to.

A \emph{RAESDC (regular, algebraic, essentially self-dual, cuspidal)
automorphic representation of $\GL_{m}(\mathbb{A}_\Q)$} can be defined as a pair $(\pi, \mu)$
consisting of a cuspidal automorphic representation $\pi$ of $\GL_{m}(\mathbb{A}_\Q)$ and a continuous character $\mu: \mathbb{A}_\Q^{\times}/ \Q^{\times} \rightarrow \C^{\times}$ such that:

\begin{enumerate}

\item (regular algebraic) $\pi$ has \emph{weight} $a= (a_i) \in \Z^n$.

\item  (essentially self-dual) $\pi \cong \pi^{\vee} \otimes (\mu \circ \mathrm{Det})$.

\end{enumerate}

Given a RAESDC automorphic representation $\pi$ as above,  there exist a number field $M\subset \C$, a finite set $S$ of rational primes, and  strictly
compatible systems of semisimple Galois representations 
\begin{equation*}\begin{aligned} \rho_\lambda (\pi): G_\Q &\rightarrow \GL_m (\overline{M}_\lambda),\\
\rho_\lambda(\mu): G_\Q &\rightarrow \overline{M}_\lambda^{\times}, \end{aligned}\end{equation*}
where $\lambda$ ranges over all finite places of~$M$
(together with fixed embeddings $M \hookrightarrow M_\lambda \hookrightarrow \overline{M}_\lambda$,
where $\overline{M}_\lambda$ is an algebraic closure of the localisation $M_\lambda$ of $M$ at $\lambda$)
such that the following properties are satisfied. Denote by $\ell$ the rational prime lying below $\lambda$.

\begin{enumerate}

\item  $\rho_\lambda(\pi) \cong \rho_\lambda(\pi)^{\vee} \otimes
\chi_{\ell}^{1-m} \rho_\lambda(\mu)$, where $\chi_{\ell}$ denotes the
$\ell$-adic cyclotomic character.

\item The representations $\rho_\lambda(\pi)$ and
$\rho_\lambda(\mu)$ are unramified outside $S \cup \{ \ell  \}$.

\item  Locally at $\ell$, the representations $\rho_\lambda(\pi)$ and
$\rho_\lambda(\mu)$ are de Rham, and if $\ell \notin S$, they are
crystalline.

\item $\rho_\lambda(\pi)$ is regular, with Hodge-Tate weights 
$\{  a_1 + (m-1), a_2 + (m-2), \ldots, a_m  \}$.

\item\label{item:4}  Fix any isomorphism $\iota:\overline{M}_\lambda\simeq \C$ compatible with the inclusion $M\subset \C$. Then
\begin{equation}\label{eq:star} \iota\WD(\rho_{\lambda}(\pi)|_{G_{\Q_p}})^{\mathrm{F-ss}} \cong \rec (\pi_p \otimes | \mathrm{Det}  |_p^{(1-m)/2}).\end{equation}
Here $\WD$ denotes the Weil-Deligne representation attached to a
representation of $G_{\Q_p}$,
$\mathrm{F}-\mathrm{ss}$ means the Frobenius semisimplification, and $\rec$ is the notation for the (unitarily normalised) Local Langlands Correspondence.
\end{enumerate}

The properties (1)--(5) above give us some information about the compatible system $(\rho_{\bullet}(\pi))$. If we want to realise groups in a given family of finite linear groups as Galois groups over $\Q$, we will need to find a suitable RAESDC automorphic representation such that the information provided by (1)--(5) allows us to ensure that the images of the corresponding residual representations $\overline{\rho}_{\lambda}(\pi)$  belong to this family. We can already make some remarks in this connection. For example, (1) implies that the image of $\rho_{\lambda}(\pi)$ lies in an orthogonal or symplectic group. (2) provides us with a strong control on the ramification of the Galois realisation that we obtain. This is a characteristic feature of this strategy of addressing the inverse Galois problem. (3) and (4) are of a technical nature, and we will not mention them in the rest of the survey (except briefly in connection to the proof of Theorem \ref{theorem:ArDiShWi}). Instead, let us expand on the last property (5). Any $\pi$ as above can be written as a certain restricted product of local components $\pi_p$, where $p$ runs through the places of $\Q$. Equation \eqref{eq:star}, with is highly involved notation, is essentially telling us that this local component $\pi_p$ determines the restriction of $\rho_{\lambda}(\pi)$ to a decomposition group $G_p\subset G_{\Q}$ at the prime $p$. As we will see in the next section, the possibility of prescribing the restriction of  $\rho_{\lambda}(\pi)$ to $G_p$ for a finite number of primes $p\not=\ell$ will be the essential ingredient for controlling the image of $\overline{\rho}_{\lambda}(\pi)$.

\section{Main statements and ingredients of proof}\label{sec:Ingredients}

In this section we state several results obtained through the
strategy described in Section \ref{sec:CompatibleSystems}, that
generalise Theorems \ref{theorem:Wiese} and \ref{theorem:DieulefaitWiese} to $2n$-dimensional representations.
The first result, due to C.~Khare, M.~Larsen and G.~Savin (cf.~\cite{KLS1}), can be encompassed in the vertical direction, as
explained in Section \ref{sec:CompatibleSystems}.

\begin{theorem}[Khare, Larsen, Savin]\label{theorem:KLS1} Fix $n, t\in \mathbb{N}$ and a prime
$\ell$. Then there exists a natural number $r$ divisible by $t$ such
that either $\PSp_{2n}(\mathbb{F}_{\ell^r})$ or
$\PGSp_{2n}(\mathbb{F}_{\ell^r})$ occurs as a Galois group over
$\mathbb{Q}$.

More precisely, there exists an irreducible Galois representation
$\rho_{\ell}:G_{\Q}\rightarrow \GSp_{2n}(\overline{\Q}_{\ell})$,
unramified outside $\ell$ and an auxiliary prime $q$, such that the
image of $\overline{\rho}^{\mathrm{proj}}_{\ell}$ is either
$\PSp_{2n}(\mathbb{F}_{\ell^r})$ or
$\PGSp_{2n}(\mathbb{F}_{\ell^r})$.
\end{theorem}

We also want to mention the following result, dealing with different families of linear groups (cf.~\cite{KLS2}).

\begin{theorem}[Khare, Larsen, Savin]\label{theorem:KLS2} Fix $t\in \mathbb{N}$ and a
prime $\ell$.

\begin{enumerate}
\item There exists an integer $r$ divisible by $t$ such that
$G_2(\mathbb{F}_{\ell^r})$ can be realised as a Galois group over
$\Q$.

\item Assume that $\ell$ is odd. There exists an integer $r$ divisible by $t$ such
that either the group $\mathrm{SO}_{2n+1}(\mathbb{F}_{\ell^r})^{\mathrm{der}}$
or  $\mathrm{SO}_{2n+1}(\mathbb{F}_{\ell^r})$ can be realised as a
Galois group over $\Q$.

\item Assume that $\ell\equiv 3, 5\pmod{8}$. There exists an integer $r$ divisible by $t$ such
that the group $\mathrm{SO}_{2n+1}(\mathbb{F}_{\ell^r})^{\mathrm{der}}$ can be
realised as a Galois group over $\Q$.

\end{enumerate}
\end{theorem}

In the horizontal direction there is the following result for symplectic groups, due to
S.~A., L.~Dieulefait, S.-W.~Shin and G.~Wiese (cf.~\cite{ArDiShWi}).

\begin{theorem}[A., Dieulefait, Shin, Wiese]\label{theorem:ArDiShWi} Fix $n, r\in \mathbb{N}$.
There exists a set of rational primes of positive density such that,
for every prime $\ell$ in this set, the group
$\PSp_{2n}(\mathbb{F}_{\ell^r})$ or $\PGSp_{2n}(\mathbb{F}_{\ell^r})$ can
be realised as a Galois group over $\Q$.

More precisely, there exists an irreducible Galois representation
$\rho_{\ell}:G_{\Q}\rightarrow \GSp_{2n}(\overline{\Q}_{\ell})$,
unramified outside $\ell$ and two auxiliary primes, such that the
image of $\overline{\rho}^{\mathrm{proj}}_{\ell}$ is either
$\PSp_{2n}(\mathbb{F}_{\ell^r})$ or $\PGSp_{2n}(\mathbb{F}_{\ell^r})$.
\end{theorem}

Note that, in \cite{DiWi}, the authors can control whether the image is $\PSL$ or $\PGL$ because they choose their modular form in such a way that it does not have any nontrivial inner twist. Currently, this has not been generalised to $n>1$.

In both results, there are essentially two different parts: on the one hand, one needs to find conditions on a compatible
system of symplectic Galois representations to ensure that the images of the residual representations corresponding to the members of the system will be huge. On the other
hand, one needs to show the existence of RAESDC automorphic representations whose compatible systems satisfy the desired
conditions. In \cite{KLS1}, the existence of appropriate automorphic representations is shown by means of Poincar\'e series, which give automorphic representations on $\SO_{2n+1}(\mathbb{A}_{\Q})$. These are transferred to $\GL_{2n}(\mathbb{A}_{\Q})$ by means of Langlands functoriality. In \cite{ArDiShWi}, the existence of the desired automorphic representations is shown by exploiting results of S.-W.~Shin on equidistribution of local components at a fixed prime in the unitary dual with respect to the Plancherel measure (cf.~\cite{Shi12}).

In the rest of the section, we will expand on the first question, namely, on conditions on symplectic compatible systems that allow some control on the images of the residual representations corresponding to the members of the system. A first property of the image that we want to ensure is irreducibility. In both \cite{KLS1} and \cite{ArDiShWi}, this is achieved by means of a tamely ramified symplectic local parameter. More precisely, fix a prime $\ell$, and let $p, q$ be auxiliary primes such that the order of $q$ modulo $p$ is exactly $2n$. Let $\Q_{q^{2n}}$ be the unique unramified extension of $\Q_q$ of degree $2n$. Using class field theory, it can be proven that there exists a character $\chi_q:G_{\Q_{q^{2n}}}\rightarrow \overline{\Q}_{\ell}^{\times}$ of order $2p$ such that (1) the restriction of $\chi_q$ to the inertia group $I_{\Q_{q^{2n}}}$ has order exactly $p$; (2) $\chi_q(\mathrm{Frob}_{q^{2n}})=-1$. Then it follows that the Galois representation $\rho_q:=\mathrm{Ind}_{G_{\Q_{q^{2n}}}}^{G_{\Q_q}}\chi_q$ is irreducible and can be conjugated to take values inside $\mathrm{Sp}_{2n}(\overline{\Q}_{\ell})$.
As a consequence, we obtain the following result:

\begin{lemma}\label{prop:irreducibility}
Let $(\rho_{\bullet})$ be a $2n$-dimensional compatible system of Galois representations of $G_{\Q}$ as in Definition \ref{defi:CompatibleSystem}. Let $p$, $q$ two primes such that the order of $q$ modulo $p$ is exactly $2n$.  Let $G_q\subset G_{\Q}$ be a decomposition group at $q$, and assume that, for all primes $\lambda$ of $L$ which do not lie above $p$ or $q$, we have
\begin{equation*}
 \mathrm{Res}^{G_\Q}_{G_{q}} \rho_{\lambda}\simeq \mathrm{Ind}_{G_{\Q_{q^{2n}}}}^{G_{\Q_q}} \chi_q,
\end{equation*}
where $\ell$ is the rational prime below $\lambda$ and $\chi_q:G_{\Q_{q^{2n}}}\rightarrow \overline{\Q}_{\ell}$ is a character as above. Then $\overline{\rho}_{\lambda}$ is irreducible.
\end{lemma}

More precisely, the image of $\overline{\rho}_{\lambda}$ contains a so-called $(2n, p)$-group (cf.~\cite{KLS1} for the definition of this notion).
Given a prime $\ell$, if one chooses the auxiliary primes $p$ and $q$ in an appropiate way, it is possible to ensure that the image of $\overline{\rho}_{\lambda}$ is huge (i.e., it contains $\Sp_{2n}(\mathbb{F}_{\ell})$). This idea appeared originally in the work of C.~Khare and J.-P.~Wintenberger on Serre's Modularity Conjecture for $n=1$, and has been exploited in \cite{Wiese2008} and \cite{KLS1}. Let us briefly sketch how it works in the case when $n=1$. Assume that we have a representation $\overline{\rho}_{\lambda}:G_{\Q}\rightarrow \GL_2(\overline{\mathbb{F}}_{\ell})$, satisfying that the restriction of $\overline{\rho}_{\lambda}$ to a decomposition group at $q$ is isomorphic to $\mathrm{Ind}_{G_{\Q_{q^{2n}}}}^{G_{\Q_q}} \overline{\chi}_q$. Consider the composition $\overline{\rho}_{\lambda}^{\mathrm{proj}}$ of $\overline{\rho}_{\lambda}$ with the projection $\GL_2(\overline{\mathbb{F}}_{\ell})\rightarrow \PGL_2(\overline{\mathbb{F}}_{\ell})$. We certainly know that the image of $\overline{\rho}_{\lambda}^{\mathrm{proj}}$ is a finite subgroup of $\PGL_2(\overline{\mathbb{F}}_{\ell})$. L.~E.~Dickson has classified all finite subgroups of $\PGL_2(\overline{\mathbb{F}}_{\ell})$ into four types of groups: a subgroup $H\subset \PGL_2(\overline{\mathbb{F}}_{\ell})$ is either (1) equal to $\PSL_2(\mathbb{F}_{\ell^r})$ or $\PGL_2(\mathbb{F}_{\ell^r})$ for some $r$; or (2) a reducible subgroup; or (3) a dihedral subgroup $D_s$ for some integer $s$ coprime to $\ell$; or (4) isomorphic to one of the alternating groups $A_4$, $A_5$ or the symmetric group $S_4$.

Since we know that the image of $\overline{\rho}_{\lambda}$  contains the subgroup $\overline{\rho}_{\lambda}(G_q)$, which is the dihedral group $D_p$, we can immediately exclude the possibilities (2) and (4) (provided $p$ is large enough so that it does not divide the cardinality of $A_5$ and $S_4$). To conclude that the image of $\overline{\rho}_{\lambda}$ is huge, we have to exclude the case that it is a dihedral group. Assume then that this is the case. Then $\overline{\rho}_{\lambda}=\mathrm{Ind}_{G_K}^{G_{\Q}} \psi$ for some quadratic field extension $K/\Q$. In addition, we know that $\mathrm{Res}_{G_q}^{G_{\Q}}\overline{\rho}_{\lambda}\simeq \mathrm{Ind}_{G_{\Q_{q^{2n}}}}^{G_{\Q_q}} \overline{\chi}_q$. Is there a way to get a contradiction? The idea is that this can be achieved, provided we choose the auxiliary primes $p$ and $q$ carefully. If this is the case, these two conditions will be rendered incompatible because of the relationship between $p$, $q$ and $\ell$. The reader interested in the details is referred to \cite{Wiese2008}. In order for this strategy to work, we must start from a prime $\ell$ and choose $p$ and $q$ accordingly. Thus, this idea is particularly well suited to address the vertical direction.

In \cite{KLS1} this idea is generalised to the $2n$-dimensional setting. The first difficulty that arises is that the classification of finite subgroups of $\GL_{2n}(\overline{\F}_{\ell})$ is much more intrincate when $n>1$. The main group-theoretical tool that is used in \cite{KLS1} is a theorem from \cite{LarsenPink2011}, which generalises a classic theorem of Jordan from characteristic zero to arbitrary characteristic. More precisely, let $m\in \mathbb{N}$ be an integer. Then there exists a constant $J(m)$ such that, for any finite subgroup $\Gamma$ of $\GL_m(\overline{\F}_{\ell})$ there exist normal subgroups $\Gamma_3\subset \Gamma_2\subset \Gamma_1\subset \Gamma$ such that the index $[\Gamma:\Gamma_1]\leq J(m)$, and such that $\Gamma_3$ is an $\ell$-group, $\Gamma_2/\Gamma_3$ is an abelian group (whose order is not divisible by $\ell$) and $\Gamma_1/\Gamma_2$ is a direct product of finite groups of Lie type in characteristic $\ell$. 

Going back to the setting of Galois representations, the main idea now is that, if $\Gamma\subset \GSp_{2n}(\overline{\F}_{\ell})$ is a finite subgroup such that there is a $(2n, p)$-group contained in all normal subgroups of $\Gamma$ of index smaller than or equal to a constant $d(n)$ which depends only on $n$ (this constant will be computed in terms of the quantity $J(2n)$ mentioned above), then it follows that $\Gamma$ must contain $\Sp_{2n}(\F_{\ell})$. Given a prime number $\ell$, by choosing the auxiliary primes $p$ and $q$ in a suitable way, one can ensure that if $\mathrm{Res}_{G_q}^{G_{\Q}} \overline{\rho}_{\lambda}\simeq \mathrm{Ind}_{G_{\Q_{q^{2n}}}}^{G_{\Q_q}} \overline{\chi}_q$, then the group $\Gamma=\mathrm{im}\overline{\rho}_{\lambda}$ satisfies that $\overline{\rho}_{\lambda}(G_q)$ is a $(2n, p)$-group contained in all normal subgroups of $\Gamma$ of index at most $d(n)$. 

Now we focus our attention on the horizontal direction. Recall that, in this setting, we are given a compatible system $(\rho_{\bullet})$, and we want that the image of the members $\overline{\rho}_{\lambda}$ are huge for as many primes $\lambda$ of $L$ as possible. In this context, the presence of a tamely ramified local parameter at an auxiliary prime $q$ will not suffice to obtain huge image. Since the prime $\ell$ is now varying, we are not allowed to choose the auxiliary primes $p$ and $q$ in terms of $\ell$. A new idea is required. 

When $n=1$, L. Dieulefait and G. Wiese construct Hecke eigenforms $f$ such that the compatible system of Galois representations $(\rho_{f, \bullet})$ attached to $f$ satisfies that, for all primes $\lambda$ of $\Q_f$, the image of $\overline{\rho}_{f, \lambda}$ is huge (cf.~\cite{DiWi}). The idea is to choose $f$ in such a way that the corresponding compatible system has \emph{two} tamely ramified parameters (at two different auxiliary primes), chosen in such a way that all possibilities for the image of $\overline{\rho}_{f, \lambda}^{\mathrm{proj}}$ given by Dickson's classification (see above) except huge image are ruled out. 

For the $2n$-dimensional case, however, we need a new ingredient. The main result in \cite{ArDiShWi} relies on a classification of finite subgroups of $\GSp_{2n}(\overline{\F}_{\ell})$ containing a transvection. More precisely, the main result in \cite{ArDiWi2} shows that, if $\Gamma\subset \GSp_{2n}(\overline{\F}_{\ell})$ is a finite subgroup which contains a (nontrivial) transvection, then either (1) $\Gamma$ is a reducible subgroup; or (2) $\Gamma$ is imprimitive; or (3) $\Gamma$ is huge. The first possibility can be ruled out by introducing a tamely ramified parameter in the compatible system $(\rho_{\bullet})$. The imprimitive case corresponds to the situation when $\rho_{\lambda}$ is induced from some field extension $K/\Q$. To rule out this case, one needs to choose the auxiliary primes $p$ and $q$ in the tamely ramified parameter in a suitable way. If the compatible system is regular (in the sense that the tame inertia weights of $\overline{\rho}_{\lambda}$ are independent of $\lambda$ and different, cf.~\cite{ArDiWi2} for a precise definition), then the second case in the classification can be ruled out, and the conclusion that the image of $\overline{\rho}_{\lambda}$ is huge can be drawn. 

The question remains whether it is possible to enforce a compatible system $(\rho_{\bullet})$ of Galois representations to satisfy, by means of a local condition, that the images of the residual representations $\overline{\rho}_{\lambda}$ contain a transvection. Recall that in Section \ref{sec:CompatibleSystems}, transvections already appeared in connection with the image of the Galois representation attached to the group of $\ell$-torsion points of an abelian variety $A$ defined over $\Q$. In this setting, to ensure that the image of $\overline{\rho}_{A, \ell}:G_{\Q}\rightarrow \GSp_{2n}(\mathbb{\F}_{\ell})$ contains a transvection, C.~Hall exploited the fact that, if $A$ has a certain type of reduction at an auxiliary prime $p_1$, then the image of the inertia group at $p_1$ under $\overline{\rho}_{A, \ell}$ already contains a transvection. In the case of $2n$-dimensional compatible systems of Galois representations, the transvection can be obtained by imposing that the restriction of $\rho_{\lambda}$ to a decomposition group at an auxiliary prime $p_1$ has a prescribed shape. Equivalently, this amounts to specifying the Weil-Deligne representation attached to the restriction of $\rho_{\lambda}$ to $G_{p_1}$. If the compatible system $(\rho_{\bullet}(\pi))$ is attached to a RAESDC automorphic representation $\pi$, this condition can be expressed in terms of $\pi$. Here it is very important that the local component $\pi_{p_1}$ of $\pi$ determines, via the Local Langlands correspondence, not only the characteristic polynomial of $\rho_{\lambda}(\mathrm{Frob}_{p_1})$ for $\lambda\nmid p_1$, but the whole restriction $\rho_{\lambda}(\pi)\vert_{G_{p_1}}$. Moreover, one has to take care that the transvection in the image of $\rho_{\lambda}(\pi)$ does not become trivial under reduction modulo $\lambda$. In \cite{ArDiShWi}, the authors ensure that, for a density one set of rational primes $\ell$ and for every $\lambda\vert \ell$, the transvection is preserved after reduction modulo $\lambda$. The main tool they use is  a level lowering result from \cite{BLGGT}, which they apply over infinitely many quadratic CM number fields.

Up to this point, we have sketched the main ideas in \cite{Wiese2008}, \cite{KLS1} and \cite{DiWi}, \cite{ArDiShWi} to prove the existence of compatible systems of Galois representations $(\rho_{\bullet})$ such that the images of the residual representations $\overline{\rho}_{\lambda}$ are huge, i.e., containing $\Sp_{2n}(\F_{\ell})$. For the applications to the inverse Galois problem, we need a certain control of the largest exponent $r$ such that  $\Sp_{2n}(\F_{\ell^r})$ is contained in the image of $\overline{\rho}_{\lambda}$. We already remarked in Section \ref{sec:Classical} that, in the case of Galois representations attached to a Hecke eigenform $f$, this is linked to the knowledge of the coefficient field $\Q_f$, which proves to be a difficult task. However, even though it may be difficult to determine precisely what the coefficient field $L$ of the compatible system is, it is possible to ensure that it contains a large subfield. In fact, the tamely ramified parameter at the prime $q$ provides already a lower bound on the size of $r$. For the applications in the horizontal direction, one exploits that if $L/\Q$ contains a cyclic subextension $K/\Q$ of degree $r$, then there exists a positive density set of primes $\ell$ such that, at some prime $\lambda$ of $L$ above $\ell$, the extension $L/\Q$ has the desired residue degree $r$.

\bibliographystyle{amsplain}

\end{document}